\documentclass[10pt]{article}
\usepackage[english]{babel}
\usepackage{latexsym}
\usepackage{color}
\usepackage{amsmath}
\usepackage[latin1]{inputenc}
\usepackage{amsfonts}
\usepackage{amssymb}
\usepackage{mathrsfs}
\usepackage{eucal}
\usepackage{exscale}
\usepackage{makeidx}
\usepackage{makeidx}
\usepackage{graphicx}
\usepackage[T1]{fontenc}
\usepackage{amsmath,amssymb}
\usepackage[all]{xy}
\newcommand{\qed}{{$\square$}\medbreak}
\newtheorem{thm}{Theorem}
\newtheorem{rem}{Remark}
\newtheorem{df}{Definition}

\newtheorem{cor}{Corollary}
\newtheorem{lem}{Lemma}
\newtheorem{fait}{Fact}
\newtheorem{exemples}{Exemples}
\newtheorem{exemple}{Exemple}

\newcommand{\preuve}{\indent {\it Proof.}\hspace{4mm}}

\title{Positive Robinson theories and h-maximal models}
\author{Mohammed Belkasmi}

\begin{document}
\maketitle


\section*{\small{Introduction}}

In this paper we continue the  exploration of  the classes
of positively closed and h-maximal model of an h-inductive 
theory in the context of positive logic.\\
 In the section 2  we give 
a concrete 
description of the class of h-maximal models of an 
h-inductive theory and theirs companion theories. The section 3
 is concerned to the study of  the positive Robinson and
 locally positive  Robinson theories and their connexion with 
 the properties of the class of h-maximal models of the  
 companion theories, and their connexion with the 
 property  of elimination 
 of quantifiers. 
Before dealing with the topics mentioned above
 we  give in section 1 a brief 
introduction to the  positive model theory.

\section{Positive model theory}

The positive logic is a continuation of the line of research on 
universal theories initiated by Abraham Robinson, 
based on the study of the notions of inductive theories, existentially closed models, model-complete theories through the
notions of embedding, existential formula. The systematic treatment of the positive model theory  has been
undertaken by Ben Yaacov and Poizat  in \cite{begnacpoizat}.\\
In short consists of non-use of negation in building
of formulas.

Let $L$ be a first order language.
The positive formulas
are  expressed as: $\exists\, \bar x\,\psi(\bar x,\bar y)$, where 
$\psi$ is a   formula, the variables  
$\bar y$ are said to be free.

A sentence is a formula without free variables. 
 A sentence is said to be h-inductive (resp. f-inductive)
 if it is a finite conjunction of sentences 
of the form: $$\forall\bar x\ \ \exists
\bar y\psi(\bar x, \bar y)\rightarrow\ \ 
\exists\bar z\varphi(\bar x,\bar z)$$
(resp. $\forall\bar x\ \ \alpha(\bar x)
\rightarrow\ \beta(\bar x)$)
where $\psi, \varphi$ (resp. $\alpha, \beta$) are 
 quantifier-free positive formulas.

The h-universal sentences represent a  special case of
h-inductive sentences, they are the sentences that can be written
as  negation of a positive sentence.

Given two L-structures $A$ and  $B$ be 
over an arbitrary language $L$. 
 A mapping $f$ from $A$ into $B$ is 
a homomorphism if for every $\bar a\in A$ and for every 
positive atomic 
formula $\phi$; $$A\models\phi(\bar a)\ \Rightarrow 
B\models\phi(\bar f(a)).$$ A structure  
$B$ is said to be  a continuation of a structure $A$ if and only if there exists  a homomorphism from $A$ into $B$.\\
 A homomorphism $f$ is an  embedding if and only if  for every 
 $\bar a\in A$ the tuples  $\bar a$ and $f(\bar a)$ satisfy the same atomic
 formulas.\\
A homomorphism $f$ is  an immersion if and only if for 
every $\bar a\in A$ and for every positive 
formula $\varphi$; $A\models\varphi(\bar a)$ 
if and only if $B\models\varphi(\bar f(a))$. We 
say that  $A$ is immersed in $B$ if there exist 
an immersion from $A$ into $B$.

A class of L-structures is said to be h-inductive if it is closed
with respect to inductive limits of homomorphisms. 
In \cite{begnacpoizat} it is shown 
 that the class of models
of an $h$-inductive theory is h-inductive and
 the class of models of an arbitrary theory $T$ 
 is h-inductive if $T$ is axiomatized by an h-inductive theory.\\

\subsection{positively  closed structures}

\begin{df}
A member $M$ of a class $\Gamma$ of $L$-structures
is said to be positively closed (pc from
now on) in $\Gamma$,
if every homomorphism from $M$ into a member  
of $\Gamma$ is an immersion.
\end{df}

\begin{fait}[{\cite[Theorem 1]{begnacpoizat}}]\label{pecconti}
Every member of an h-inductive class of $L$-structures
has a positively  closed continuation in the same class.
\end{fait}
The h-inductivity of the class $\Gamma$  is a necessary 
condition of the existence of 
pc structure. In this case 
The class of pc members of $\Gamma$ forms an h-inductive and h-cofinal
subclass of $\Gamma$.

Let $\Gamma$ be an h-inductive class of $L-structure$. We denote 
by $\Pi(\Gamma)$ the class of positively closed member of $\Gamma$. 
If $\Gamma$ is the class of models of an h-inductive theory
$T$, we use the notation $\Pi(T)$.

\begin{df}
Two $h$-inductive theories over a language $L$ are said 
to be companions if 
they have the same pc models.
\end{df}

Note that every  $h$-inductive
theory $T$ admits:
\begin{itemize}
\item  A maximal companion theory denoted $T_k(T)$,  called
the Kaiser's hull theory of $T$. By definition $T_k(T)$
is the set of $h$-inductive sentences satisfied by the pc models of
 $T$.
\item A minimal companion theory denoted $T_u(T)$, it is the
set of  h-universal sentences 
true in the pc models of $T$.
\end{itemize}

Let $L$ be a first order language and $M$ be a $L$-structure. 
\begin{itemize}
\item we denote by $T_i(M)$ (resp. $T_v(M)$) the set of 
h-inductive (resp. of h-universal) sentences satisfied by $M$ 
in the language obtained from $L$ by adding the elements of $M$ as constants.
\item we denote by $T_i^*(M)$ (resp. $T_v^*(M)$) the set of 
h-inductive (resp. of h-universal) $L$-sentences satisfied by $M$.
\end{itemize}
Note that for every L-structure $M$ we have;
$$T_i^{\star}(M)\subset T_k(T_i^{\star}(M)) ,\ \  \ T_v^{\star}(M) \subset T_u(T_v^{\star}(M)).$$
In the language obtained from $L$ by adding the elements of $M$ as 
constants. We have;
$$T_k(T_i(M))= T_i(M),\ \ \ \ T_u(T_v(M))= T_v(M).$$
If $A$ a pc model of an h-inductive theory $T$. We obtain; 
$$T_i^{\star}(A)= T_k(T_i^{\star}(A)) ,\ \ \ T_v^{\star}(A)= T_u(T_v^{\star}(A)).$$ 
In this case we use the notation $T_k^\star(A)$ instead of $T_i^\star (A)$.

\begin{df}
Let $T$ be an h-inductive theory.
\begin{itemize}
 \item $T$ is said to be model-complete 
if  every model of $T$ is a pc model of $T$.
\item We say that $T$ has a model-companion whenever $T_k(T)$ is 
model-complete.
\item An $n$-type is a maximal set of positive formulas
in $n$ variables that is consistent with $T$. We denote by $S_n(T)$  the space
of $n$-types of a theory $T$. 
\end{itemize}
 \end{df}
Let $M$ be a L-structure and $\bar m$ a tuple of $M$ . We denote by  
$tp_M(\bar n)$ the set of positive formulas satisfied by $\bar m$
in $M$.
\begin{fait}
 $M$ is pc model of $T$  if and only if, for every $\bar a\in A$,
the set of positive formulas satisfied by $\bar a$
is a type of $T$.
\end{fait}
For every positive formula $\phi$, we denote by   
$Ctr_T(\phi)$ the  set of  positive formulas $\psi$ such that 
$T\vdash\neg\exists x(\phi(\bar x)\wedge\psi(\bar x))$.

Let  $A$ be a  pc model of $T$. Let
$\bar a\in A$ such that $A\nvDash\phi(\bar a)$
where $\phi$ is a positive formula. By the maximality of 
$tp_A(\bar a)$, there is 
  a positive formula $\psi\in Ctr_T(\phi)$ such that $A\models\psi(\bar a)$.\\
This property is in fact the inner characteristic of these  subclass of 
models of $T$. 
We have the following fact.
\begin{fait}
 $A$ is pc model of $T$  if and only if for every 
$\bar a\in A$, and for every positive formula $\varphi$;
if $A\nvDash\varphi(\bar a)$  there exists a positive formula $\psi$
such that $A\models\psi(\bar a)$ and $\psi\in Ctr_T(\phi)$.
\end{fait}

Consider a pc model  $A$  of $T$ and 
$\bar a\in A$. We denote by $tp(\bar a)$ (resp. $tpqf(\bar a)$)
the type of $\bar a$ in $A$ (resp. the set of quantifier-free
positive formulas satisfied by $\bar a$ in $A$).

One defines on $S_n(T)$ the  topology generated by the following 
 basis of closed sets: 
\[
F_\varphi\ =\ \{\ p\in S_n(T)\ |\ p\vdash\varphi\ \}\ .
\]
where $\varphi$ ranges over the 
set of positive formulas.\\
Note that for every $n$, The space of positive types $S_n(T)$ is compact
but generally  is not  Hausdorff.

\begin{df}
Let $T$ be an h-inductive theory and $\varphi$ a positive 
formula;
\begin{itemize}
\item $\varphi$ is said to be $T$-complemented  
if and only if there is 
a positive formula $\psi\in Ctr_T(\varphi)$ such that;
$$T\vdash\forall\bar x\ \ (\varphi(\bar x)\vee\psi(\bar x)).$$ 
The formula $\psi$ is called the $T$-complement of $\varphi$.
\item Let $\Gamma$ be a subset of $Ctr_T(\varphi)$. We say that 
$Ctr_T(\varphi)$ is logically equivalent to $\Gamma$ modulo 
$T$ and we writ 
$Ctr_T(\varphi)\approx_T \Gamma$;
if and only if for every $\psi\in Ctr_T(\varphi)$ there is 
$\phi\in\Gamma$ such that 
$$T\vdash\forall\bar x\ (\psi(\bar x)\rightarrow\phi(\bar x)).$$
\end{itemize}
\end{df}
\begin{rem}
Let $T$ be an h-inductive theory.
\begin{itemize}
\item A formula $\varphi$ is 
$T_k(T)$-complemented if and only if there is a positive formula 
$\psi$ such that  $Ctr_T(\varphi)\approx_{T_k(T)}\psi$.
\item  The class of pc models of $T$ is elementary
if and only if, for every positive formula $\psi$, $Ctr_T(\psi)$
is logically equivalent modulo $T_k(T)$ to a  positive 
formula.
\end{itemize}
 
\end{rem}

\begin{exemples}\label{exemple}
\begin{enumerate} 
\item Let $L$ be the relational language formed a binary relation $S$.
Consider the following h-inductive theory:
$$T=\{\neg\exists\,xy\, (S(x,y)\wedge S(y,x)), 
\forall xyz\ ((S(x, z)\wedge S(y, z))\rightarrow x=y)\}.$$
The model of $T$
formed by the p-cycles where $p=4$ or 
 $p$ is a prime number greater-than or equal to $3$ is the unique 
 pc model of $T$.\\
 Let $T'$ be the theory obtained from $T$
By adding  the  h-universal sentence
$$\neg\exists x_1x_{2}x_{3}x_{4}\  
((\bigwedge_{i=1}^{3} S(x_i, x_{i+1}))\wedge S(x_4, x_1)).$$

The structure formed by the  $p$-cycles where $p$
ranges over the set of prime numbers  greater-than or equal to $3$ 
is the unique pc model of $T'$.
\item Let $L$ and $T'$ be the language and the theory given in 
the example above. Let $n$ be an integer greater-than 3. 
Consider $T_n$ the h-inductive theory obtained from $T'$ by adding 
the following set of h-inductive sentences
$$\{\forall x_1\cdots x_m\, 
((\bigwedge_{i=1}^{m} S(x_i, x_{i+1}))\wedge S(x_m, x_1))
\longrightarrow \bigvee_{i\neq j} x_i=x_j\ |\  m>n\}.$$
The structure formed by the p-cycles
where  $p$ is a prime number  less than  $n$
is the unique pc model of 
$T_n$. Thereby $T_n$ has a model-companion.
\item Let $T_{ag}$ be the h-inductive theory of 
abelian groups in the language 
$L=\{.,\, ^{-1}, e\}$.
In the positive logic $T_{ag}$ has a model-companion.
The trivial group $\{e\}$ 
is the unique pc model of  $T_{ag}$. 
However,
 in the context of first order logic 
the class of existentially closed abelian groups is
the class of divisible abelian groups which contain for
each  prime $p$  an infinite number 
elements of order $p$ (theorem 2.4 \cite{module}).
\end{enumerate}
\end{exemples}
 
To extend the discussion began on the last example.
Consider the language $L^\star$  obtained from the language $L$ of 
the theory $T_{ag}$ by  
adding a constant $a$. Let  $T_{ag}^+$ be 
the h-inductive theory 
$T_{ag}, \{\neg a=e\}$.\\ 
 Let 
$(G, g)$ be a pc model of $T_{ag}^+$ where $g$ is the interpretation 
of the constant $a$ in $G$, we have the following properties
\begin{enumerate}
\item The constant $g$ belongs to every non trivial  
subgroup of $G$.
Indeed let $N$ be a non trivial  subgroup of $G$
 and $\pi$ the L-homomorphism $G\rightarrow G/N$.  
Suppose  that $\pi$ is a $L^\star$-homomorphism. Then $\pi$ 
is an immersion. Consequently $N=\{e\}$. Thereby $\pi$
can not be a $L^\star$-homomorphism,  so $\pi(g)=e$.\\
  The constant $g$ belongs to the  intersection 
of all subgroups of $G$. Thereby for every $x\in G$ there is 
$k\in \mathbb{Z}$ such that $g=x^k$.

\item  $G$ cannot admit distinct subgroups of order $p$ and  $q$ respectively,
where $p$ and $q$ are prime to each other.
Because if not, the order of $g$ will be a common divisor of 
$p$ and $q$.
\item $G$ cannot be the direct sum of some of its subgroups; 
because the constant must belong to the intersection of all 
subgroups.
\end{enumerate} 
\begin{lem}
The  pc models of $T_{ag}^+$ are the groups 
$\{(Z(p),z_p)$, where $p$ is a prime number,
 $z_p\in \{Z(p)-1\}$, and 
$Z(p)$ is the group of all complex $p^n$-th roots of unity.
\end{lem}
\preuve 
Let $(G, g)$ be a pc model of  $T_{ag}^+$. We distinguish two cases:
\begin{itemize}
\item $o(g)$ (the order of $g$ in $G$) is finite.
 In this case $o(g)$ is a prime
number. Indeed, if not  we can find  a subgroup of $H=<g>$ 
 that  does not contain the constant $g$.
\item $o(g)$ is infinite. This case  cannot take place because
 we can find a subgroup of $H=<g>$  which does not 
contain the constant $g$.
\end{itemize}
Therefore, if $(G, g)$ is a pc model of $T_{ag}^+$ there is $p$ a prime number
such that $(G, g)$ is the group of all complex $p^n$-th roots of unity.

 \subsection{Amalgamation property}
The notion of amalgamation in positive logic provides  a useful means
for intuition and motivation. One of these facts is the characterization of the Hausdorff property by the amalgamation 
property given in \cite{begnacpoizat}. For more expositions of 
these facts see \cite{ana, begnacpoizat}. 
\begin{df}
Let $\Gamma$ be a class of L-structures. An element $A$ of $\Gamma$ 
is said to be  an amalgamation basis of $\Gamma$
 if and only if, for every $B, C$ in $\Gamma$, and $f, g$ 
homomorphisms 
from $A$ respectively into $B$ and $C$, there exist 
$D\in\Gamma$ and $f', g'$ homomorphisms 
 such that
the following diagram commutes:
\[
\xymatrix{
    A \ar[r]^{f} \ar[d]_{g} & {B} \ar[d]^{g'} \\
    C \ar[r]_{f'} & {D}
  }
\]
 We say that $\Gamma$ has the amalgamation 
property if every element of $\Gamma$ is an amalgamation basis of $\Gamma$.
\end{df} 

Note that under certain conditions, each structure can  
benefit of the property of being an amalgamation basis. 
On other words in every class
of L-structures, we can always find  universal amalgamations.
 The useful following fact provides an example of these 
 universal amalgamations.
\begin{fait}[{\cite[lemma 4]{ana}}]\label{asyamalg}
Let $A, B, C$ be L-structures such that; $A$ is immersed in $B$
and continued in $C$ by a homomorphism $f$. Then there is 
$D$ a model of $T_k(C)$ such that the following diagram commutes.
$$\xymatrix{
    A \ar[r]^{i_m} \ar[d]_{f} & B \ar[d]^{g} \\
    C \ar[r]_{i_m} & {D}
  }$$
  Where $i_m$ in the diagram are immersions and $g$ a homomorphism.
\end{fait}
One of the most important property of the class
 of pc models of 
an h-inductive theory is the amalgamation property 
(theorem 9\cite{begnacpoizat}). As a simple application 
of the amalgamation property  we have the following
lemma.
\begin{lem}\label{singleton}
Let $T$ be an h-inductive theory such that the class  of 
pc model (resp. of amalgamation bases) is closed under product.
 Then $T_k(T)$ 
 has only one pc model, this pc model has only one point. 
\end{lem} 
\preuve 
Let $A, B$ be two pc models of $T$. 
Since $A\times B$ is a pc model (resp. an amalgamation basis),
we obtain the following commutative diagram:
$$\xymatrix{
    A\times B \ar[r]^{pr_A} \ar[d]_{pr_B} & A \ar[d]^{f} \\
    B \ar[r]_{g} & {C}
  }$$
where $C$ is a pc model of $T$,  $f$ and $ g$ are immersions. Thus
for all $a\in A$ and $b\in B$ we have $f(a)=g(b)$, and so
 $tp(a)=tp(b)$. Consequently every constant mapping  from 
$A$ into $A$ is an immersion. Thereby $A=\{a\}$.

\subsection{Complete theories} 
\begin{df}
An h-inductive theory $T$ is said to be complete or has 
joint continuation property
if any two of its models can be simultaneous continued into a
third one.
\end{df}
 \begin{rem}\ 
 \begin{itemize}
 \item let $A, B$ be two pc models of $T$.
$T_k^\star(A)=T_k^\star(B)$ if and only if $A$ and $B$ have
 the same continuation. 
 \item an h-inductive theory is complete if and only if  its 
 pc models have the same h-inductive theory.
 \end{itemize}

 \end{rem}
\begin{lem}\label{peccomplete}
Let $A$ be a pc model of $T$ then $T_k^\star(A)$ is a complete theory.
\end{lem}
\preuve 
Let $B, C$ be models of $T_k^\star(A)$. Firstly we  show that\\
$\{T_k^\star(A), Diag^+(A), Diag^+(B)\}$ is consistent, 
then we conclude that\\ 
$\{T_k^\star(A), Diag^+(B), Diag^+(C)\}$ is consistent.\\
Since for every $\varphi(\bar a)\in Diag^+(A)$ we have 
$\exists\bar x\varphi(\bar x)\in T_k^\star(A)$, 
and $B\vdash T_k^\star(A)$, then 
by compactness  we obtain the  consistency of 
$\{T_k^\star(A), Diag^+(A), Diag^+(B)\}$.\\
Now, let $B^\star$ be  a model of 
$\{T_k^\star(A), Diag^+(A), Diag^+(B)\}$ 
and $C^\star$  a model of $\{T_k^\star(A), Diag^+(A), Diag^+(C)\}$. 
Since 
the class of pc models of $T$ has the amalgamation property,
$B^\star$ and $C^\star$ are  models of $T$. We obtain 
the following commutative diagram:
$$\xymatrix{
B \ar[r]& B^\star  \ar[rd]^{f_1}&\\
A  \ar[ru]^{i_1} \ar[rd]_{i_2}&& D\\
C \ar[r]& C^\star  \ar[ru]_{f_2}& \\
     }$$
where $i_1, i_2$ are immersions, $f_1, f_2$ are  homomorphisms
and $D$ is a model
of $T$ that  can be assumed a pc model of $T$. 
Thus,  $f_1\circ i_1$ is an 
immersion, and   $D$ is a model 
of $T_k^\star(A)$ in which $B$ and $C$ are immersed. 
Consequently 
$\{T_k^\star(A), Diag^+(B), Diag^+(C)\}$ is consistent, and  
$T_k^\star(A)$ is complete theory.
\begin{cor}
Let $A$ be a pc model of $T$. Every pc model of $T_k^\star(A)$ is a pc model 
of $T$.
\end{cor}
\preuve
Let $C$ be a pc model of $T_k^\star(A)$ and $B$ a pc model of $T$ 
such that $C$  
is continued in $B$ by a homomorphism $f$. By  the lemma 
\ref{peccomplete}  and the fact\ref{asyamalg} we obtain the following commutative diagram:
$$\xymatrix{
A \ar[r]^{i'}& B'  \ar[rd]^{g}&&\\
C  \ar[ru]^{i} \ar[rd]_{f}&& D\ar[r]^{h}& D'\\
    & B  \ar[ru]_{i_m}& \\
     }$$
where $B'$ is a model of $T_k^\star(A)$ in which $A$ and $C$ are immersed
(lemma \ref{peccomplete}).  
$D\vdash T_i(B)$  
(fact\ref{asyamalg}). 
$D'$ is a pc model of $T$ in which $D$
 is continued by a homomorphism $h$.\\
Now, since $D'$ and $A$ are pc models
of $T$, then $h\circ g\circ i'$ is an immersion and 
 $D'\vdash T_k^\star(A)$.
 By the fact that $C$ is a pc model of $T_k^\star(A)$ we deduce 
that  $h\circ i_m\circ f$  is an immersion. Consequently $f$ is an immersion and 
     $C$ is a pc model of $T$.

 \begin{lem}
Let $T_1$ and $T_2$ be two h-inductive theories such that $T_1$ is 
a complete theory and there exists a common pc model $A$ of $T_1$ and 
$T_2$, then every pc model of $T_1$ is a pc model of $T_2$.
\end{lem} 
\preuve 
Let $B$ be a pc model of $T_1$. Since $A$ is a common 
pc model of $T_1$ and $T_2$, and $T$ is a complete theory, then
$$T_k(T_2)\subset T_k^\star(A)= T_k^\star(B).$$
 This implies that $B$ is a model of $T_2$.\\
 On the other hand, since $T_1$ is complete there exist 
 $D$ a pc model of $T_1$ in which $A$ and $B$ are immersed, and  
so $D$ is a model of $T_2$. Let $C$ be a pc model of 
$T_2$ in which $D$ is continued by an homomorphism $f$ as shown 
in the following diagram
\[
\xymatrix{
    &B \ar[d]^{i_m} &   \\
    A \ar[r]_{i_m} & D\ar[r]_{f}& C
  }
\]
where $i_m$ denotes immersions. Given that  $A$ and $C$ are  pc models 
of $T_2$ we obtain
$T_k^\star (A)= T_k^\star (C)$. Consequently   $C\models T_1$ and $f$ is an 
immersion, so $D$ and $B$ are pc models of $T_2$.
\section{H-maximal models}\label{hmax}
In \cite{almaz} Kungozhin  introduced the notion of h-maximal
 model in the context of studying
the elementarity of the classes of pc modeles and h-maximal models 
of finitely universal theories. In this section ????

\begin{df}
Let $T$ be an h-inductive theory. A model $A$ of $T$ is said
to be h-maximal if every homomorphism from $A$ into a
model of $T$ is an embedding.
\end{df}
Note that  
the class of h-maximal models of an h-inductive theory $T$ forms 
an inductive class  and every model of $T$ is continued
 in a h-maximal 
 model of $T$.
\begin{exemples}
\begin{enumerate}
\item Let $T$ and $T'$ be the theories defined in [1, example\ref{exemple}].
 The class of 
h-maximal models of $T$ (resp. $T'$) is the class 
of substructures of the pc model 
of $T$ (resp. $T'$).
\item The class of h-maximal models of the theory $T_n$ given in 
[2, example \ref{exemple}] is the class of substructures 
of the pc model of $T_n$. This implies that  the class of h-maximal models 
of $T_n$ is not elementary.\\

\item The class of h-maximal models of $T_{ag}^+$ is the class 
of pc models of $T_{ag}^+$.\\

\item Consider $T_{g}^+=T_{ag}^+-\{\forall xy\ xy=yx\}$ 
the theory  of groups in the language 
$L^\star$.
The h-maximal models of $T_{g}^+$  are the groups  whose 
non trivial normal subgroups contain the constant  of the language
$L^\star$. Indeed, since the $L^\star$-homomorphisms are 
The  homomorphisms $f$ of  
groups such that $f(a)\neq 0$ where $a$ is the interpretation constant of
 $L^\star$. Thus if $(G, a)$ is a  h-maximal model of $T_{g}^+$
  and $N$ a non trivial normal subgroup of $G$ such that $a\notin N$.
  The canonical mapping from $G$ into $G/N$ is a
   $L^\star$-homomorphism but  not an embedding. Thereby
   $(G, a)$ can not be h-maximal model.
\end{enumerate}
\end{exemples}
\begin{rem}
 If $\Gamma_1$ and  $\Gamma_2$ are h-inductive
 classes that have the same h-maximal models. Then they are 
 companion theories.
 \end{rem}
Let $T$ be an h-inductive theory and $T_m(T)$ the set of 
h-inductive sentences  satisfied in each h-maximal model 
of $T$. Given that the class of  pc models of $T$ is a subclass
of the class of h-maximal models of $T$, then
 $T\subseteq T_m(T)\subseteq T_k(T)$. So
  $T$ and $T_m(T)$ are companion theories. 

\begin{df}
Let $T$ be an h-inductive theory and $\Gamma$ 
the set of sentences of the form 
$\forall\bar x\,\phi(\bar x)\rightarrow\psi(\bar x)$ 
 satisfied in every h-maximal model of $T$ and
such that $\phi$ and $\psi$ are quantifier-free positive 
formulas. We denote by $T_f(T)$ the h-inductive 
theory $T, \Gamma$. we have 
$$T\subseteq T_f(T)\subseteq T_m(T).$$
\end{df}
\begin{lem}
The h-inductive theories 
  $T, T_f(T)$ and $T_m(T)$ have
 the same class of h-maximal models. 
\end{lem} 
\preuve
Let $A$ be a h-maximal model of $T$, $B$ a model of $T_f(T)$
and $C$ a model of $T_m(T)$. Let $f$ (resp. $g$) be 
 a homomorphism from $A$ into 
$B$ (resp. into $C$).
 Since  
$B$ and $C$ are models of $T$,
and  $A$ is a model of both theories
$T_m(T)$ and $T_f(T)$ 
then $f$ and $g$ are 
embeddings, and  $A$ is a h-maximal of $T_f(T)$ and $T_m(T)$.

Let $A$ be a h-maximal 
model of $T_m(T)$. Given that $A$ is a model
of the theories $T_f(T)$ and $T$, there are $B$ a model of $T$ and 
$C$ a model $T_f(T)$ such that $A$ is continued 
in $B$ by a homomorphism $f$  and continued in $C$ by a 
homomorphism $g$. Since 
$B$ and $C$ are also models of $T$, there exist $B'$ 
and $C'$ h-maximal models of $T$ such that $B$ 
is continued in $B'$ by a homomorphism $f'$, $C$
is continued  in 
$C'$ by a homomorphism $g'$. Now given that $B'$ and 
$C'$ are models of $T_m(T)$ then $f'\circ f$ and $g'\circ g$ are 
embeddings. Thereby $g$ and $f$ are embeddings. Consequently
$A$ is a h-maximal of $T$ and $T_f(T)$.\\
By the same way we show that every h-maximal of $T_f(T)$ is 
a h-maximal of $T$. Therefore the theories $T, T_f(T)$ and $T_m(T)$
have the same class of h-maximal models.
\begin{rem}
Consider $T$  an h-inductive theory. We denote by $\Sigma_T$ the class 
of h-maximal model of $T$. We have 
$\Sigma_{T_k(T)}\subseteq \Sigma_{T}=\Sigma_{T_f(T)}=\Sigma_{T_m(T)}\subseteq\Sigma_{T_u}$.

\end{rem}

\begin{lem}\label{carach-max}
A model $A$ of $T$ is  h-maximal model  if and only if 
for every quantifier-free positive formula  $\varphi$ and  a tuple
 $\bar a\in A$
 such that 
$A\nvDash \varphi(\bar a)$, there is
  $\psi(\bar x)$ a 
positive formula $\psi\in Ctr_T(\varphi)$
such that $A\models\psi(\bar a)$.\\
\end{lem}
\preuve
Let $A$ be a h-maximal model of $T$, $\bar a\in A$ and  $\varphi$ 
a quantifier-free positive formula such that 
$A\nvDash \varphi(\bar a)$.
 Since every homomorphism from $A$
into a model of $T$ is an embedding then the  set of h-inductive sentences  
$\{T, Diag^+(A), \varphi(\bar a) \}$ is inconsistent. Thus
by compactness there exists $\phi(\bar a, \bar b)\in Diag^+(A)$
such that 
$T\vdash \neg\exists\bar x(\varphi(\bar x)\wedge\psi(\bar x))$ 
where $\psi$
is the positive formula $\exists \bar y\phi(\bar x,\bar y)$.

Conversely, let $A$ be a model of $T$ such that for every quantifier-free positive 
formula $\varphi$ and $\bar a\in A$, if $A\nvDash \varphi(\bar a)$
then there is   $\psi(\bar x)$ a 
positive formula 
such that $A\models\psi(\bar a)$ and 
$T\vdash\neg\exists\bar x(\varphi(\bar x)\wedge\psi(\bar x))$. It is obvious  
that every homomorphism from $A$ into a model of
 $T$ is an embedding, then $A$ is a h-maximal model of $T$.
\begin{cor}
If $A$  is immersed in a h-maximal model of $T$ 
then $A\in \Sigma_T$.
\end{cor}
\preuve Since $A$ is immersed in a model of $T$ then $A\vdash T$. The fact
that $A$ is h-maximal results of the lemma \ref{carach-max}.\qed

\begin{thm}\label{hmaxetpec}
 $\Sigma_T$
 is elementary class if and only if for every quantifier-free positive 
formula $\varphi$, there is a positive formula $\psi$ such that 
$$Ctr_{T_m(T)}(\varphi)\approx_{T_m(T)}\{\psi\}.$$
\end{thm}
\preuve
Suppose that $\Sigma_T$ is elementary and axiomatized by $T_m(T)$.
Assume  the existence of  a quantifier-free positive formula
 $\varphi$ such that 
$Ctr_{T_m(T)}(\varphi)$ is not equivalent modulo $T_m(T)$
to any positive formula.  
By compactness, there is $B$ 
  a model of $T_m(T)$ and $\bar b\in B$ 
such that   $B\nvDash\varphi(\bar b)$, and
for every positive formula $\psi\in Ctr_{T_m(T)}(\varphi)$ we have 
  $B\nvDash\psi(\bar b)$, which contradicts 
the lemma \ref{carach-max}.

For the reverse direction, suppose that for every quantifier-free positive 
formula $\varphi$, there is a positive formula 
$\psi\in Ctr_{T_m(T)}(\varphi)$ such that 
$T_m(T)\vdash\forall x\, \varphi(x)\vee\psi(x)$. Let 
$A$ be a model of $T_m(T)$. By the lemma \ref{carach-max} it is clear
that $A$ is a h-maximal model of $T$.
\begin{cor}
If $\Sigma_{T} $ is elementary then $\Sigma_{T_k(T)} $
is elementary and axiomatized by $T_k(T)$.
\end{cor}
\preuve Suppose that $\Sigma_{T} $  is 
axiomatised by $T_m(T)$.  by the theorem (\ref{hmaxetpec}), 
for every quantifier-free positive 
formula $\varphi$ there is a positive formula 
$\psi\in Ctr_{T_m(T)}(\varphi)$
 such that 
$T_m(T)\vdash\forall x\, \varphi(x)\vee\psi(x)$. Given that
$T_k(T)\supseteq T_m(T)$,
then every model of $T_k(T)$ is a h-maximal model of $T_k(T)$.\qed

\begin{lem}
If $\Sigma_{T_k(T)} $ is elementary then  it is axiomatized by 
$T_k(T)$.
\end{lem}
\preuve Suppose that $\Sigma_{T_k(T)} $ is axiomatized 
by an h-inductive theory $T^\star$. Then $T^\star$ and 
$T_k(T)$ are companion theories. Given that 
 $T_k(T)$ is the maximal companion of 
$T$ we obtain $T_k(T) \sim T^\star$.
\section{Positive Robinson and locally positive  Robinson theories}
In \cite{rob} Hrushovski defined Robinson theories to be  the universal theory 
that admits the quantifier separation.   The quantifier-free
types are the main object of the study 
of Robinson theories. In our context we adopt this property
 to define the notion of positive Robinson theories and 
 locally positive Robinson theories.
\begin{df}\label{posrobdefn} 
An $h$-inductive theory $T$ is said to be positive Robinson theory
(in short. pR theory) if it satisfies the following condition:\\
For any  pc models $A$ and $B$ of $T$,  $\bar a\in A$ and
$\bar b\in B$. If  
$tpqf(\bar a)\subseteq tpqf(\bar b)$ then $tp(\bar a)=tp(\bar b)$.
(where $tpqf(\bar a)$ is the set of quantifier-free positive formulas satisfied
by $\bar a$ in $A$).\\
An h-inductive theory $T$ is said to be a 
locally positive Robinson theory (in short. lpR theory) if
the following conditions is satisfied for any   
 pc model $A$ of $T$. 
$$\forall \bar a, \bar b\in A;\ \ tpqf(\bar a)\subseteq tpqf(\bar b)
\Rightarrow tp(\bar a)=tp(\bar b)$$
\end{df}
Given that the property of being a pR theory or 
a lpR  theory concerns the class of pc models. we have
the following remarks.
\begin{rem}\label{comprob}
\begin{itemize}
\item  $T$ is a  
pR theory (resp. lpR theory)  provided that each companion 
theory  of $T$ is a pR theory (resp. lpR theory).
\item  If $T$ is a pR  theory then  $T$ is a lpR theory.
\end{itemize}
\end{rem}

\begin{fait}\label{caracparres}[Lemma 8, \cite{ana}]
An $h$-inductive theory $T$ is a pR theory if and only if
for every positive formula $\varphi$, $Ctr_T(\varphi)$ is equivalent modulo
$T_k(T)$ to a set of quantifier-free positive formulas.
\end{fait}
\begin{rem}
Let $T$ be a  pR theory. If $A\in \Sigma_T$
and $B$ a model of $T$ which is embedded in $A$. Then 
$B\in\Sigma_T$.
\end{rem}
\begin{thm}
 $T$ is lpR h-inductive theory if and only if for every pc model 
$A$ of $T$, and $\varphi$ a positive 
formula  we have the following property:\\
for every  tuple $\bar a\in A$,  if $A\nvDash\varphi(\bar a)$ 
then there  exists  $\psi$ a free 
positive formula such that, 
$A\models\psi(\bar a)$ and 
$T_k^\star(A)\vdash\neg\exists\bar x\ \ (\varphi(\bar x)\wedge\psi(\bar x))$.  
\end{thm}
\preuve Suppose that $T$ is a lpR theory.
Let $A$ be a pc model of $T$, $\bar a\in A$ and $\varphi$  a
positive formula such that, $A\nvDash\varphi(\bar a)$. We will show
that $T^\star =\{T_k^\star(A), tpqf_A(\bar a), \varphi(\bar a)\}$ 
is  inconsistent, where 
$tpqf_A(\bar a)$ is the set of quantifier-free positive formulas 
satisfied by $\bar a$ in the pc model $A$.\\
Suppose that $T^\star$ is consistent. Let 
$B$ a model of $T^\star$ in the language 
$L^\star=\{L, \{\bar a\}\}$. We claim 
that $\{T, Diag^+(A), Diag^+(B)\}$ is
 consistent. Indeed if not, 
 by compactness there exist $\psi(\bar a, \bar x)\in Diag^+(A)$  
 a quantifier-free positive formula
 such that 
 $\{T, Diag^+(B), \psi(\bar a, \bar x)\}$ is inconsistent.
 Given that $B\vdash \{T, Diag^+(B)\}$ then 
 $B\vdash\neg\exists\bar x\ \psi(\bar a, \bar x)$. On 
 the other hand, since $B\vdash T_k^\star(A)$ and 
 $A\models\exists\bar x\ \psi(\bar a, \bar x)$ we obtain 
 $B\models\exists\bar x\ \psi(\bar a, \bar x)$, contradiction.
 Thereby  $\{T, Diag^+(A), Diag^+(B)\}$ is consistent.
 Let $C$ be a model of $\{T, Diag^+(A), Diag^+(B)\}$, 
  so $A$ and $B$ are continued in $C$. Let $\bar b$   the interpretation 
 of $\bar a\in B$ in $C$, and $\bar a$
 the interpretation  of $\bar a\in A$ in $C$.
  Given that  every model 
 of $T$ is continued in some pc model of $T$, we can take
 $C$ a pc model of $T$.
 \\
Considering that $A$ is a pc model of $T$, it is immersed in $C$,  
thereby we obtain $tpqf_A(\bar a)=tpqf_C(\bar a)$. 
Since $B\vdash T^\star$, we have
$$tpqf_A(\bar a)\subseteq tpqf_B(\bar a)\subseteq tpqf_C(\bar b).$$
Since $T$ is lpR theory we have
 $tp(\bar a)=tp(\bar b)$.\\
Given that  $A\nvDash\varphi(\bar a)$ and $B\models\varphi(\bar b)$ 
 then $C\nvDash\varphi(\bar a)$ and $C\models\varphi(\bar b)$.
Contradiction. Therefore $T^\star$ is inconsistent, by 
compactness there exists $\psi(\bar x)\in tpqf_A(\bar a)$ such 
that 
$T_k^\star(A)\vdash\neg\exists\bar x\ (\psi(\bar x)
\wedge\varphi(\bar x))$.

For the reverse direction, suppose that for every pc model $A$ of $T$,
$\bar a\in A$, and $\varphi$ a positive formula
we have the following property:\\ if 
 $A\nvDash\varphi(\bar a)$ then there exists a quantifier-free
positive formula $\psi$ such that 
$A\models\psi(\bar a)$, and 
$T_k^\star(A)\vdash\neg\exists\bar x\ \ 
(\varphi(\bar x)\wedge\psi(\bar x))$.
Let $\bar a$ and  $\bar b$ be  tuples of $A$ such that 
$tpqf_A(\bar a)=tpqf_A(\bar b)$. Assume the existence of  
a positive formula $\varphi$ such that $A\models\varphi(\bar a)$ and 
$A\nvDash\varphi(\bar b)$.  By hypothesis there is  a 
quantifier-free positive formula $\psi$ such that 
$A\models\psi(\bar b)$ and 
$T_k^\star(A)\vdash\neg\exists\bar x\ (\varphi(\bar x)\wedge\psi(\bar x))$.
Thus $\psi\in tpqf(\bar b)$, but we have $\psi\notin tpqf(\bar a)$. 
Contradiction. Thereby $T$ is lpR theory.
\begin{lem}
An h-inductive theory $T$ is lpR if and only if for every
pc model $A$ of $T$, $T_k^\star(A)$ is a pR theory.
\end{lem} 
 \preuve
 Suppose that $T$ is a lpR theory. Let $A$ be a pc model 
 of $T$. Let  $B$ and $C$ be two pc models of $T_k^\star(A)$.  
Consider $\bar b\in B$ and
 $\bar c\in C$ such that $tpqf_B(\bar b)\subseteq tpqf_C(\bar c)$.
 By   
 the lemma \ref{peccomplete}, the sets 
  $S_1=\{T_k^\star(A), Diag^+(B), Diag^+(A)\} $ and 
  $S_2=\{T_k^\star(A), Diag^+(C), Diag^+(A)\}$ are consistent. 
From  the fact \ref{asyamalg}  we obtain the following diagram:
 $$\xymatrix{
B \ar[r]^{e_1}& B^\star  \ar[rd]^{f_1}&&&&\\
A  \ar[ru]^{i_1} \ar[rd]_{i_2}&& D\ar[r]^{j_1}&D_1\ar[r]^{j_2}&D_2\\
C \ar[r]^{e_2}& C^\star  \ar[ru]_{f_2}&&&& \\
     }$$
 where $B^\star\models S_1$ and $C^\star\models S_2$ that can be 
 taken pc models of $T_k^\star(A)$. $e_1, e_2, i_1, i_2, f_1, f_2$
 are immersions,
 $D$  a model of $\{T_k^\star(A), Diag^+(B), Diag^+(C)\}$ and 
 $D_1$ a pc model of $T$ in which $D$ is continued.
 Given that $A$  and $D_1$  are pc models of $T$ 
  then $j_1\circ f_1\circ i_1$ is an immersion, which implies that 
 $D_1$ is a model of $T_k^\star(A)$. Let $D_2$ be
 a pc model of $T_k^\star(A)$ in which $D_1$ is continued.\\
 Since $j_1\circ f_1\circ e_1$ and $j_1\circ f_2\circ e_2$ 
 are immersions, we have;  
 $$tpqf_{D_1}(\bar b)= tpqf_{B}(\bar b)\subseteq tpqf_{C}(\bar c)
 = tpqf_{D_1}(\bar c)$$
 Given that $T$ is a lpR theory and $D_1$ is a pc model of 
 $T$ we obtain $tp_{D_1}(\bar b)=tp_{D_1}(\bar c)$. On the other hand, as 
 $D_2$ is a pc model of $T_k^\star(A)$ and $D_1$ is immersed in $D_2$,
 then $D_1$ is a pc model of $T_k^\star(A)$, thereby 
 $tp(\bar b)=tp(\bar c)$, wich implies that $T_k^\star(A)$ is a pR theory. \\
 The other direction of the proof is obvious.
 
 \begin{fait}\label{caracrobin}[Lemma, \cite{ana}]
Let $T$ be a  pR theory. We have the following properties:\\
- Every model of $T$ that embeds in a pc model of 
$T$ is $h$-maximal model of $T$.\\
- The $h$-maximal models of $T$ have the amalgamation property.\\
Moreover, if $T$ is $h$-universal the two conditions above
imply that $T$ is a  pR theory.
\end{fait}

From  the fact \ref{caracrobin} and the definition of pR
theories   we obtain the following slightly modified version of 
 the fact\ref{caracrobin}. 
\begin{lem}\label{caracrobamalgamation}
An h-inductive theory $T$ is pR if and only if the 
class of substructures of the h-maximal models of $T$ has the amalgamation 
property.
\end{lem}
\preuve Suppose that $T$ is a pR theory, then $T_u(T)$ the 
h-universal companion  of $T$ is also a pR theory.
Given that the class of substructures of the h-maximal models of 
$T$ is  $\Sigma_{T_u(T)}$ the class of h-maximal model of $T_u(T)$,
by the fact \ref{caracrobin},
$\Sigma_{T_u(T)}$ has the amalgamation property.

For the other direction of the proof, assume that the class 
$\Sigma_{T_u(T)}$ has the amalgamation property. 
By the fact \ref{caracrobin}, $T_u(T)$ and thus $T$ are pR theories.

\begin{cor}
 If the class of pc model of $T$ is closed under substructures 
 then  $T$ is a pR theory. 
\end{cor}
\preuve  Since
every pc model is a h-maximal  model and the class of pc  models has
the amalgamation property, the proof follows from 
the lemma \ref{caracrobamalgamation}
\begin{cor}
$T_{ag}^+$ is a pR theory.
\end{cor}
\begin{exemples}
\begin{enumerate}
\item 
Let $T$ be an h-inductive theory  that satisfies the following 
property:\\ For every pc model $A$ of $T$ and
for every positive 
formula $\varphi$,  there exists  
 a family of free-positive formulas $\{\phi_i\ |\ i\in I \}$
 such that:
\begin{itemize}
\item for every $i\in I$ we have
 $T_k^\star(A)\vdash\ \forall\bar x\ \ \phi_i(\bar x)\ \rightarrow
 \varphi(\bar x)$.
 \item $\forall \bar a\in A$, if $A\models\varphi(\bar a)$ then 
 there is $i\in I$ such that $A\models\phi_i(\bar a)$.
\end{itemize}
We claim that $T$ is a lpR theory, and 
for every pc model $A$ of $T$ the class of h-maximal models of 
$T_k^\star(A)$ is the class of pc model of $T_k^\star(A)$.\\
Indeed, let $A$ be a pc model of $T$,  $\bar a$ and 
$\bar b$ be  tuples
from $A$ such that $tpqf_A(\bar a)\subseteq tpqf_A(\bar b)$. 
Let $\varphi$ be a positive formula such that
 $A\models\varphi(\bar a)$, by hypothesis  there is $\phi$ a quantifier-free positive 
 formula such that  $T_k^\star(A)\vdash\forall\bar x\ \phi(\bar x) \rightarrow
 \varphi(\bar x)$, and $A\models\phi(\bar a)$. Given that 
$tpqf_A(\bar a)\subseteq tpqf_A(\bar b)$ and $\phi\in tpqf_A(\bar a)$,
then $A\models\phi(\bar b)$ and thereby $A\models\varphi(\bar b)$. 
Consequently $tp_A(\bar a)\subseteq tp_A(\bar b)$, as 
$A$ is a pc, by maximality of types we obtain 
$tp_A(\bar a)=tp_A(\bar b)$.\\
Now, let $B$ be a h-maximal of $T_k^\star(A)$ and  $C$ a pc model of 
$T_k^\star(A)$ in which $B$ is embedded. Let $\varphi$ a positive 
formula such that $C\models\varphi(\bar b)$ where $\bar b\in B$, 
then there is $\phi$ a quantifier-free positive formula such that: 
$$T_k^\star(A)\equiv T_k^\star(C)\vdash\forall\bar x\ \phi(\bar x) \rightarrow
 \varphi(\bar x)\ \ \text{and}\ \ C\models\phi(\bar b).$$ Since 
 $B$ is embedded in $C$ and $\phi$ is a quantifier-free formula then
$B\models\phi(\bar b)$, so $B\models\varphi(\bar b)$. 
Consequently $B$ is immersed in $C$ which implies that  $B$ is a pc model 
of $T_k^\star(A)$.

 Theories whose pc are finite provide a concrete example of  theories with the above property.
\item Let $T$ be an h-inductive theory such that for every 
positive formula $\varphi$ there is $\{\phi_i|\ i\in I\}$ a 
family of quantifier-free positive formula such that\\
\begin{itemize}
\item $T_k(T)\vdash\ \forall\bar x\ \phi_i(\bar x)\rightarrow\varphi(\bar x)$.
\item For every pc model $A$ of $T$ and $\bar a$ a tuple of $A$. If 
$A\models\varphi(\bar a)$ then $A\models\phi_i(\bar a)$ for same 
$i\in I$.
\end{itemize} 
Then $T$ is a pR theory, and every h-maximal 
model of $T_k(T)$ is a pc model of $T$.
\item Let $L$ be the language 
formed by the function symbol $f$ or arity 1. Let 
$T$ be the h-universal theory 
$\{ \neg\exists x\ f(x)=x\}$. In (\cite{almaz}, example 3) 
it is shown that the unique  pc model of $T$ is the model 
formed by the $p$-cycle (cycle of length $p$) where $p$ 
 runs through the set of prime numbers. The h-maximal models 
of $T$ are the substructures of the pc model of $T$.\\
It is obvious that $T$ is a pR  theory.
 \end{enumerate} 
\end{exemples}
 
\begin{thm}\label{th4}
Let $T$ be a pR theory with a model-companion.  
Then every positive formula is equivalent  modulo $T_k(T)$ to a 
quantifier-free positive  formula.
\end{thm}
\preuve
Since $T$ is  pR  theory with a model-companion  then for each 
positive formula $\varphi$ there is $\phi$ a quantifier-free positive 
formula such that 
$$T_k(T)\vdash\neg\exists\varphi(\bar x)\wedge\phi(\bar x),\ \ 
T_k(T)\vdash\forall\bar x\,\varphi(\bar x)\vee\phi(\bar x).$$
 We repeat the same reasoning for the quantifier-free positive formula 
 $\phi$ and we obtain a quantifier-free positive formula 
 $\psi$  such that,  
$$T_k(T)\vdash\neg\exists\phi(\bar x)\wedge\psi(\bar x), 
T_k(T)\vdash\forall\bar x\,\phi(\bar x)\vee\psi(\bar x)$$
 which implies that
$T_k(T)\vdash\forall\bar x\varphi(\bar x)\leftrightarrow\psi(\bar x)$.

\begin{cor}
Let $T$ be a lpR theory such that 
for every pc model $A$ of $T$, the theory $T_k^\star(A)$ 
has a model-companion.
Then  every positive formula  is equivalent 
modulo $T_k^\star(A)$ to a quantifier-free positive formula.
\end{cor}
\begin{exemple}
For every pc model $Z_p$ of $T_{ag}^+$.
Every positive formula in the language of the theory $T_{ag}^+$
is equivalent modulo  $T_k^\star(Z_p)$ to a quantifier-free positive formula.
\end{exemple}

\bibliographystyle{plain}

Mohammed  Belkasmi\\
Math  Department\\
College of Science\\
Qassim University
\end{document}